\documentclass[11pt,a4paper,reqno]{amsart}
\usepackage{hyperref,graphics,amsmath}
\newtheorem{theorem}{Theorem}[section]
\newtheorem{remark}[theorem]{Remark}

\theoremstyle{theorem}

\newtheorem{lemma}[theorem]{Lemma}
\numberwithin{equation}{section} \makeatletter
\@namedef{subjclassname@2010}{ 2010 Mathematics Subject
Classification} \makeatother
\title[On the polar derivative of a polynomial]{On the polar derivative of a polynomial}
\author{N. A. Rather, S. H. Ahangar and Suhail Gulzar}
 \address{Department of Mathematics \\
    University of Kashmir \\
   Srinagar, Hazratbal 190006
   \\ India}
 \email{dr.narather@gmail.com}
 \email{sgmattoo@gmail.com}
 \email{ahangarsajad@gmail.com}
\date{}
\begin{document}
\maketitle{}
\begin{abstract}
  Let $P(z)$ be a polynomial of degree $n$ having no zero in $|z|<k$ where $k\geq 1,$ then for every real or complex number $\alpha$ with $|\alpha|\geq 1$ it is known 
\begin{equation*}
\underset{|z|=1}{\max}|D_\alpha P(z)|\leq n\left(\dfrac{|\alpha|+k}{1+k}\right)\underset{|z|=1}{\max}|P(z)|,
\end{equation*} 
where $D_\alpha P(z)=nP(z)+(\alpha-z)P^{\prime}(z)$ denote the polar derivative of the polynomial $P(z)$ of degree $n$ with respect to a point $\alpha\in\mathbb{C}.$ In this paper, by a simple method, a refinement of above inequality and other related results are obtained. 
\end{abstract}
\begin{center}
\section{\textbf{ Introduction and statement of results}}
\end{center}
\hspace{10mm} If $P(z) $ is a polynomial of degree $n,$ then concerning the estimate of $ |P^{\prime}(z)| $ on the unit disk $|z|=1,$ we have
\begin{equation}\label{e1}
\underset{\left|z\right|=1}{\max}\left|P^{\prime}(z)\right|\leq  n\, \underset{\left|z\right|=1}{\max}\left|P(z)\right|.
\end{equation} 
 Inequality \eqref{e1} is an immediate consequence of Bernstein's inequality on the derivative of a trigonometric polynomial (for reference, see \cite{mm}, \cite{rs} or \cite{asc}). Equality in \eqref{e1} holds for  $P(z)=a z^n,$ $a\neq 0.$ \\
\indent If we restrict ourselves to the class of polynomials having no zero in $|z|<1$, then inequality \eqref{e1} can be replaced by
\begin{equation}\label{e2}
\underset{\left|z\right|=1}{\max}\left|P^{\prime}(z)\right|\leq        \frac{n}{2}\underset{\left|z\right|=1}{\max}\left|P(z)\right|.
\end{equation}
 Inequality \eqref{e2} was conjectured by Erd\"{o}s and later verified by Lax \cite{el}.The result is sharp and equality holds for $P(z)=\alpha z^n+\beta,$  $|\alpha|= |\beta|\neq 0.$\\
\indent As an extension of \eqref{e2}, Malik \cite{m} proved that if $P(z)$ is a polynomial of degree $n$ which does not vanish in $|z|<k$ where $k\geq 1,$ then
\begin{equation}\label{e3}
\underset{\left|z\right|=1}{\max}\left|P^{\prime}(z)\right|\leq        \frac{n}{1+k}\underset{\left|z\right|=1}{\max}\left|P(z)\right|.
\end{equation}
The result is best possible and equality in \eqref{e3} holds for $P(z)=(z+k)^n.$\\
\indent Let $D_\alpha P(z)$ denote the polar derivative of the polynomial $P(z)$ of degree $n$ with respect to the point $\alpha\in\mathbb{C},$ then
$$D_\alpha P(z)=nP(z)+(\alpha-z)P^{\prime}(z). $$
The polynomial $D_\alpha P(z)$ is of degree at most $n-1$ and it generalizes the ordinary derivative in the sense that
$$\underset{\alpha\rightarrow\infty}{\lim}\dfrac{D_\alpha P(z)}{\alpha}=P^{\prime}(z).$$

\indent  A. Aziz \cite{a88} extended inequality \eqref{e3} to the polar derivative and proved that if $P(z)$ is a polynomial of degree $n$ having no zero in $|z|<k$ where $k\geq 1,$ then for $\alpha\in\mathbb{C}$ with $|\alpha|\geq 1$ 
\begin{equation}\label{e5}
\underset{|z|=1}{\max}|D_\alpha P(z)|\leq n\left(\dfrac{|\alpha|+k}{1+k}\right)\underset{|z|=1}{\max}|P(z)|.
\end{equation} 
The result is best possible and equality in \eqref{e5} holds for the polynomial $P(z)=(z+1)^n.$\\ 
\indent The bound in \eqref{e5} depends only upon the modulus of the zero of smallest modulus and not on the moduli of other zeros. It is interest to obtain a bound which depends upon the location of all the zeros rather than just on the location of the zero of smallest modulus.\\

\indent In this paper, by a simple method, we first present the following result which is a refinement of inequality  \eqref{e5}. \\

\begin{theorem}\label{t1}
Let $P(z)=a_n\prod_{\nu=1}^{n}(z-z_\nu)$ be a polynomial of degree $n.$ If $|z_\nu|\geq k_\nu\geq 1$ where $1\leq \nu\leq n,$ then for $\alpha\in\mathbb{C}$ with $|\alpha|\geq 1,$
\begin{equation}\label{t1e}
\underset{|z|=1}{\max}|D_\alpha P(z)|\leq n\left(\dfrac{|\alpha|+t_0}{1+t_0}\right)\underset{|z|=1}{\max}|P(z)|,
\end{equation}
where 
\begin{align}\label{t0}
t_0=
\begin{cases}
1+\dfrac{n}{\sum_{\nu=1}^{n}\frac{1}{k_\nu-1}}& \text{if }\quad k_\nu> 1 \quad\text{for all}\quad\nu, \,1\leq \nu\leq n\\
1              & \text{if}\quad k_\nu=1 \quad\text{for some}\quad \nu,\, 1\leq \nu\leq n.
\end{cases}
\end{align}
\end{theorem}

\begin{remark}
\textnormal{
 If $k_\nu\geq k,$ $k\geq 1$ for $1\leq \nu\leq n,$ then $t_0\geq k$ which implies 
$$ \dfrac{|\alpha|+t_0}{1+t_0}\leq \dfrac{|\alpha|+k}{1+k}\quad\textnormal{for}\quad|\alpha|\geq 1. $$
This shows \eqref{t1e} is refinement of inequality \eqref{e5}.
} 
\end{remark}
\begin{remark}
\textnormal{If we divide the two sides of inequality \eqref{t1e}  by $|\alpha|$ and letting $|\alpha|\rightarrow \infty,$ we get a result due to Govil et. al. \cite{gl}.}
\end{remark}
Next, as an application of Theorem \ref{t1}, we present the following result.\\

\begin{theorem}\label{t2}
Let $P(z)=a_n\prod_{\nu=1}^{n}(z-z_\nu)$ be a polynomial of degree $n$ with $P(0)\neq 0.$ If $|z_\nu|\leq k_\nu\leq 1,$ $1\leq \nu\leq n,$ then for $\delta\in\mathbb{C}$ with $|\delta|\leq 1,$
\begin{equation}\label{t2e}
\underset{|z|=1}{\max}|D_{\delta} P(z)|\leq n\left(\dfrac{1+|\delta|\,s_0}{1+s_0}\right)\underset{|z|=1}{\max}|P(z)|,
\end{equation}
where 
\begin{align}\label{v0}
s_0=
\begin{cases}
1+\dfrac{n}{\sum_{\nu=1}^{n}\frac{k_\nu}{1-k_\nu}}& \text{if }\quad k_\nu< 1 \quad\text{for all}\quad\nu, \,1\leq \nu\leq n\\
1              & \text{if}\quad k_\nu=1 \quad\text{for some}\quad \nu,\, 1\leq \nu\leq n.
\end{cases}
\end{align}
\end{theorem} 

\begin{remark}
\textnormal{If $k_\nu\leq k\leq 1$ for $1\leq \nu\leq n,$ then $1/k\leq s_0$ which implies 
$$ \dfrac{1+|\delta|\,s_0}{1+s_0}\leq \dfrac{|\delta|+k}{1+k}\quad\textnormal{for}\quad|\delta|\leq 1. $$
Therefore, it follows that if $P(z)$ is a polynomial of degree $n$ having all its zeros in $|z|<k$ where $k\leq 1,$ then for $\delta\in\mathbb{C}$ with $|\delta|\leq 1,$
\begin{equation}
\underset{|z|=1}{\max}|D_\delta P(z)|\leq n\left(\dfrac{|\delta|+k}{1+k}\right)\underset{|z|=1}{\max}|P(z)|.
\end{equation}
The result is sharp.
} 
\end{remark}

\begin{center}
\section{\textbf{Lemmas}}
\end{center}

For the proof of these theorems, we need the following lemmas. The first lemma is due to  Gardner and Govil \cite{gg93}.\\

\begin{lemma}\label{l1}
Let $P(z)=a_n\prod_{\nu=1}^{n}(z-z_\nu)$ be a polynomial of degree $n.$ If $|z_\nu|\geq k_\nu\geq 1,$ $1\leq \nu\leq n,$ then for $|z|=1,$
\begin{equation}\label{l1e}
|Q^{\prime}(z)|\geq t_0|P^{\prime}(z)|,
\end{equation}
where  $Q(z)=z^n\overline{P(1/\overline{z})}$ and $t_0$ is given by \eqref{t0}.
\end{lemma}

\begin{lemma}\label{l2}
Let $P(z)$ be the polynomial of degree $n$ and $Q(z)=z^n\overline{P(1/\overline{z})}$ then for $|z|=1,$
\begin{equation}\label{l2e}
|P^{\prime}(z)|+|Q^{\prime}(z)|\leq n\,\underset{|z|=1}{\max} |P(z)|.
\end{equation}  
\end{lemma}
\noindent This is a special case of a result due to Govil and Rahman \cite{gr69}.

\begin{center}
\section{\textbf{Proof of Theorems}}
\end{center}

\begin{proof}[\textnormal {\textbf{Proof of Theorem \ref{t1}}}]
Let $ Q(z)=z^n\overline{P(1/\overline{z})}, $ then for $|z|=1$, it can be easily verified that 
\begin{align}\label{p1e1}
|P^{\prime}(z)|=|nQ(z)-zQ^{\prime}(z)|\quad\textnormal{and}\quad |Q^{\prime}(z)|=|nP(z)-zP^{\prime}(z)|. 
\end{align}
Now, for every real or complex number $\alpha$ and $|z|=1,$ we have by using \eqref{p1e1}, 
\begin{align}\nonumber\label{p1e2}
|D_\alpha P(z)|&=|nP(z)+(\alpha-z)P^{\prime}(z)|\\\nonumber&\leq |\alpha||P^{\prime}(z)|+|nP(z)-zP^{\prime}(z)|\\&=(|\alpha|-1)|P^{\prime}(z)|+|P^{\prime}(z)|+|Q^{\prime}(z)|.
\end{align}
Multiplying the two sides of inequality \eqref{p1e2} by $t_0$ and using Lemma \ref{l1}, we obtain for $|\alpha|\geq 1,$ 
\begin{align}\label{p1e3}\nonumber
t_0 |D_\alpha P(z)|&\leq (|\alpha|-1)t_0|P^{\prime}(z)|+t_0(|P^{\prime}(z)|+|Q^{\prime}(z)|)\\&\leq (|\alpha|-1)|Q^{\prime}(z)|+t_0(|P^{\prime}(z)|+|Q^{\prime}(z)|)\quad\textnormal{for}\quad|z|=1.
\end{align}
Adding \eqref{p1e2}, \eqref{p1e3} and using Lemma \ref{l2}, we get for $|\alpha|\geq 1$ and $|z|=1,$ 
\begin{align}
(1+t_0)|D_\alpha P(z)|&\leq (|\alpha|+t_0)(|P^{\prime}(z)|+|Q^{\prime}(z)|)\\\nonumber&\leq n(|\alpha|+t_0)\underset{|z|=1}{\max}|P(z)|,
\end{align}
which gives 
\begin{align}
|D_\alpha P(z)|\leq n\left(\dfrac{|\alpha|+t_0}{1+t_0}\right)\underset{|z|=1}{\max}|P(z)|\quad\textnormal{for}\quad|z|=1.
\end{align}
This completes the proof of theorem \ref{t1}.
\end{proof}

\begin{proof}[\textnormal{\textbf{Proof of Theorem \ref{t2}}}]
Since $P(z)=a_n\prod_{\nu=1}^{n}(z-z_\nu)$ where $|z_\nu|\leq k_\nu\leq 1,$ $\nu=1,2,\cdots,n$ with $P(0)\neq 0,$  $Q(z)=z^n\overline{P(1/\overline{z})}=\overline{a_n}\prod_{\nu=1}^{n}(1-\bar{z}_\nu z)$ is a polynomial of degree $n$ with $\left|\frac{1}{z_\nu}\right|\geq \frac{1}{k_\nu}\geq 1.$ Applying Theorem \ref{t1} to the polynomial $Q(z)$ and noting that $|P(z)|=|Q(z)|$ for $|z|=1,$ we obtain for $\alpha \in\mathbb{C}$ with $|\alpha|\geq 1$ and $|z|=1,$ 
\begin{align}\label{p1}
\underset{|z|=1}{\max}|D_\alpha Q(z)|\leq n\left(\dfrac{|\alpha|+s_0}{1+s_0}\right)\underset{|z|=1}{\max}|P(z)|,
\end{align}
where 
\begin{align*}
s_0=
\begin{cases}
1+\dfrac{n}{\sum_{\nu=1}^{n}\frac{k_\nu}{1-k_\nu}}& \text{if }\quad k_\nu< 1 \quad\text{for all}\quad\nu, \,1\leq \nu\leq n\\
1              & \text{if}\quad k_\nu=1 \quad\text{for some}\quad \nu,\, 1\leq \nu\leq n.
\end{cases}
\end{align*}
Again, for $|z|=1$ so that $z\overline{z}=1,$ we have
\begin{align*}\nonumber\label{daq}
|D_\alpha Q(z)|&=|nQ(z)+(\alpha-z)Q^{\prime}(z)|\\\nonumber&=\left|nz^n\overline{P(1/\overline{z})}+(\alpha-z)\left\{nz^{n-1}\overline{P(1/\overline{z})}-z^{n-2}\overline{P^{\prime}(1/\overline{z})}\right\}\right|\\\nonumber&=\left|\alpha\left\{nz^{n-1}\overline{P(1/\overline{z})}-z^{n-2}\overline{P^{\prime}(1/\overline{z})}\right\}+z^{n-1}\overline{P^{\prime}(1/\overline{z})}\right|\\\nonumber&=\left|\alpha\left(n\overline{P(z)}-\overline{z}\overline{P^{\prime}(z)}\right)+\overline{P^{\prime}(z)}\right|\\&=|\overline{\alpha}nP(z)+(1-\overline{\alpha}z)P^{\prime}(z)|=|\overline{\alpha}||D_{1/\overline{\alpha}}P(z)|.
\end{align*}
Hence, from \eqref{p1}, we get
\begin{align}
|\alpha|\underset{|z|=1}{\max}|D_{1/\bar{\alpha}} P(z)|\leq n\left(\dfrac{|\alpha|+s_0}{1+s_0}\right)\underset{|z|=1}{\max}|P(z)|.
\end{align}
Replacing $1/\overline{\alpha}$ by $\delta,$ we obtain for every real or complex number $\delta$ with $|\delta|\leq 1,$
\begin{align}
\underset{|z|=1}{\max}|D_{\delta} P(z)|\leq n\left(\dfrac{1+|\delta|\,s_0}{1+s_0}\right)\underset{|z|=1}{\max}|P(z)|.
\end{align}
This completes the proof of theorem \ref{t2}.
\end{proof}
\textbf{Acknowledgement.}
The authors are highly grateful to the referees for careful reading and their valuable suggestions. The second and third authors are supported by Council of Scientific and Industrial Research, New Delhi.

\end{document}